\documentclass[12pt]{amsart}
\usepackage{psfig}
\date{}
\begin{document}

\title[]
{On an error in the star puzzle by Henry E. Dudeney}
\author{Alex Ravsky}
\email{oravsky@mail.ru}
\address{Department of Functional Analysis, Pidstryhach Institute for Applied Problems of Mechanics and Mathematics,
Naukova 2-b, Lviv, 79060, Ukraine}
\begin{abstract}
We found a solution of the star puzzle (a path on a chessboard from c5 to d4 in 14 straight strokes)
in 14 queen moves, which has been claimed by the author as impossible.
\end{abstract}

\maketitle \baselineskip15pt

Recently I have spent many nice hours with the old puzzle book \cite{Dud1} by Henry E. Dudeney.

The task of the star puzzle, ((X36 = AM329) according to the list \cite{Knu} by Donald E. Knuth)
published a century ago in The Strand Magazine \cite{Dud2}, is to construct a path,
consisting of 14 straight strokes, on the following field of stars from one light
star to the other such that all the stars lay on the path.

\begin{figure}[h]
\hskip330pt
\psfig{file=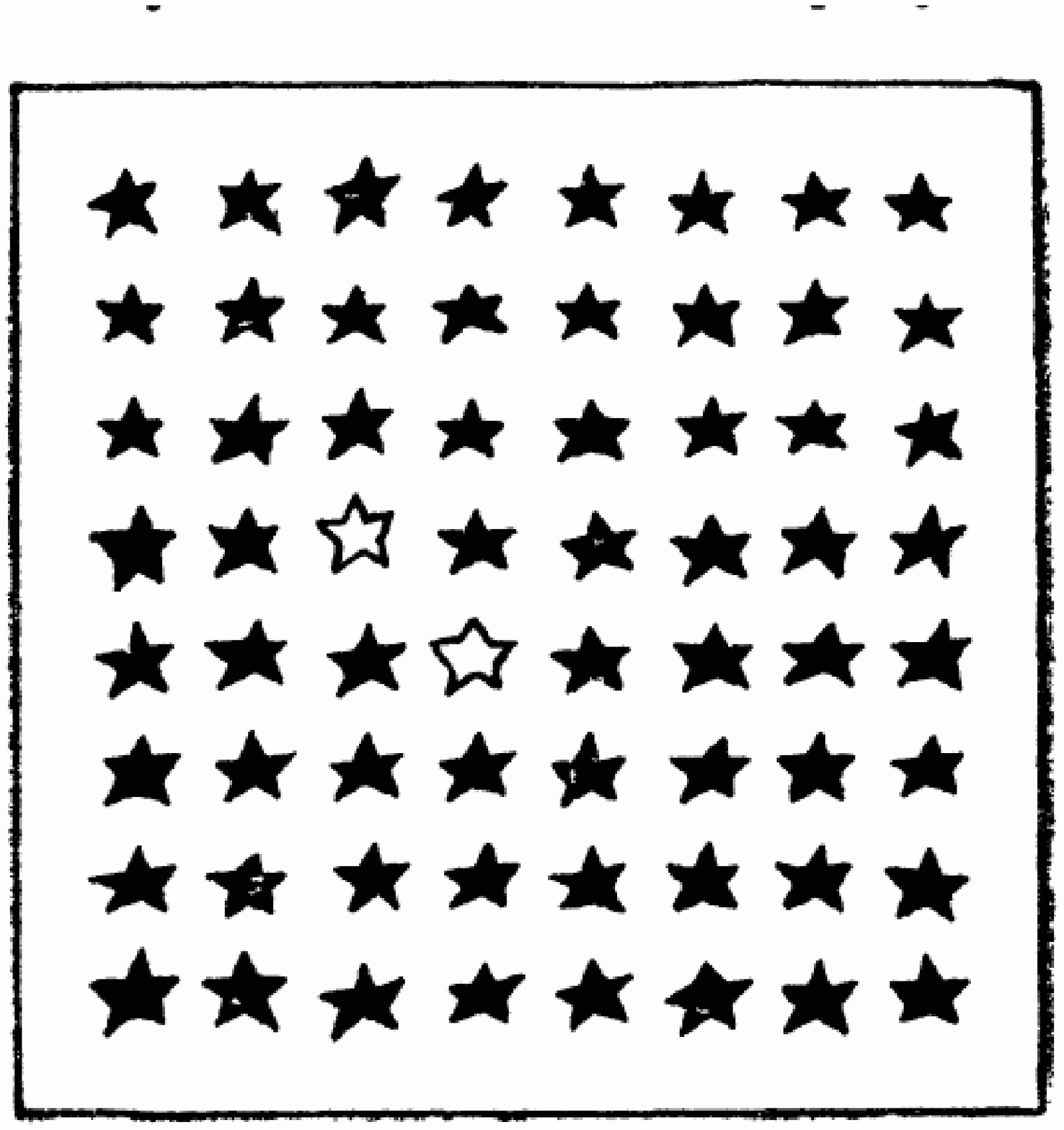,width=50mm}
\end{figure}

Author have claimed that it is impossible to do in 14 queen moves
and have proposed the following solution:

\begin{figure}[h]
\hskip330pt
\psfig{file=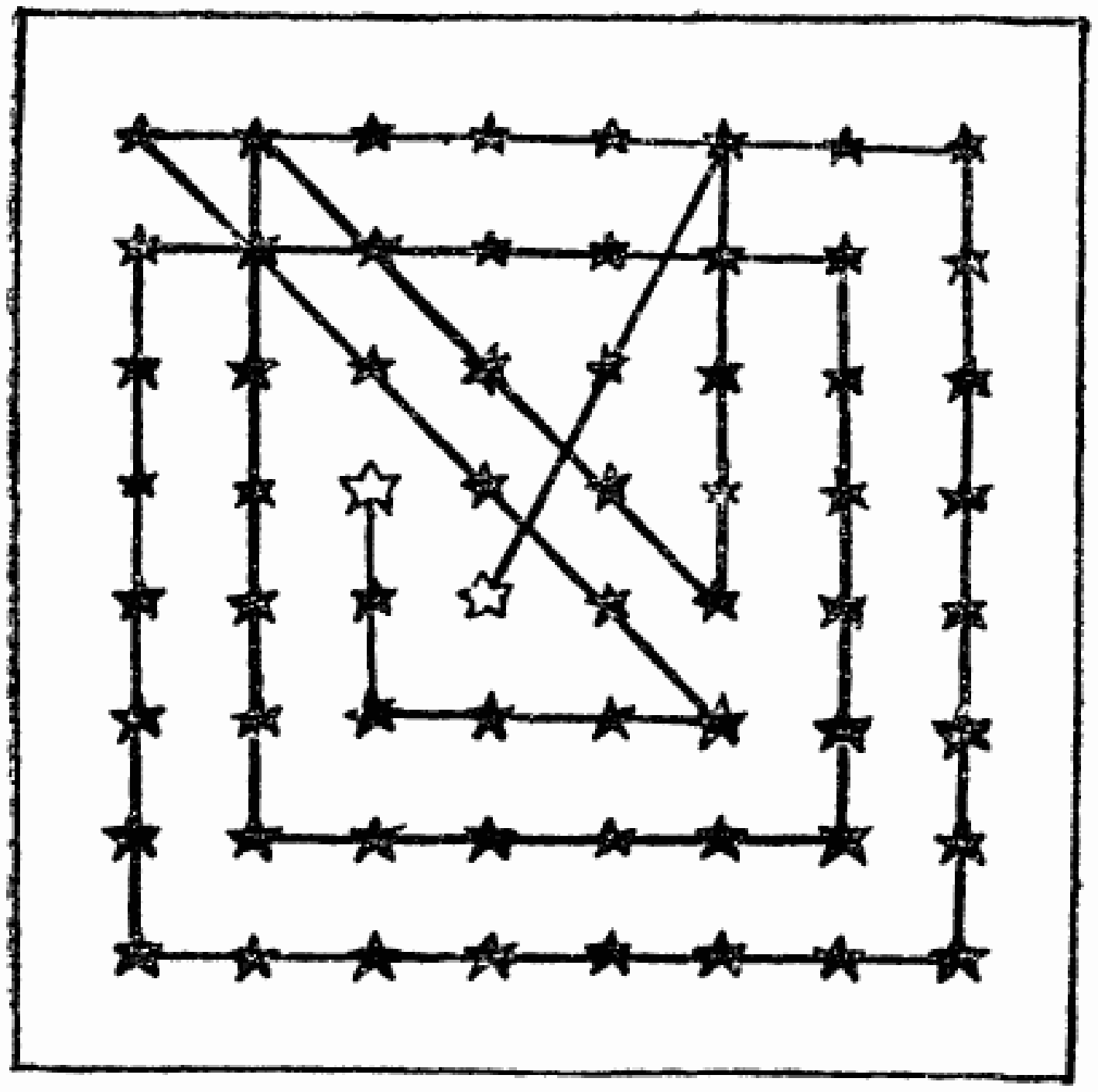,width=50mm}
\end{figure}

Surprisingly, I have found the following solution in 14 queen moves, searching it by hand,
as in the good old days:
\pagebreak

\begin{figure*}[h]
\hskip330pt
\psfig{file=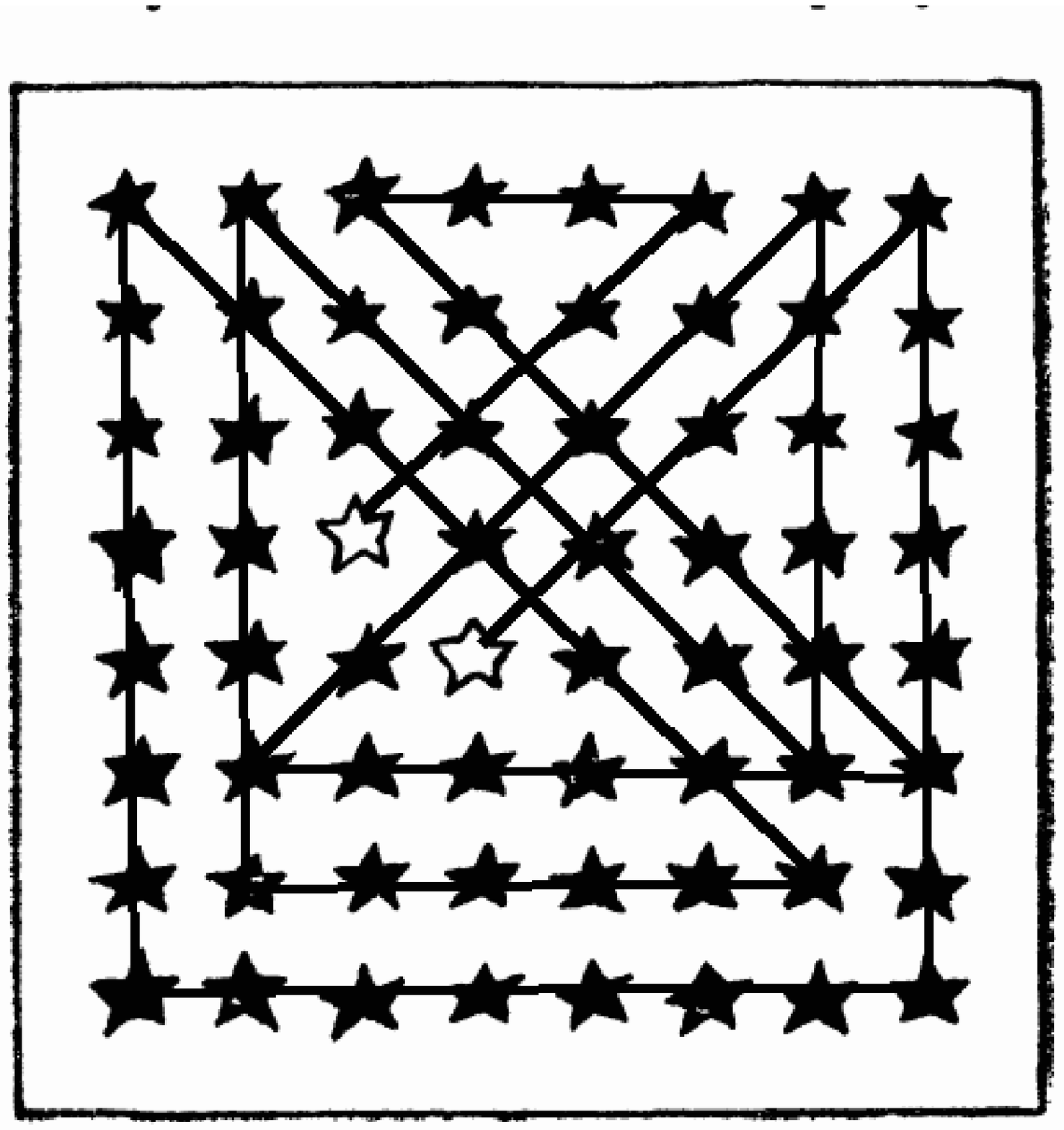,width=50mm}
\end{figure*}
\centerline{c5-f8-c8-h3-b3-g8-g3-b8-b2-g2-a8-a1-h1-h8-d4}

{}
\end{document}